\documentclass[12pt]{article}

\usepackage{graphicx}
\usepackage{enumerate}
\usepackage[vlined,ruled]{algorithm2e}
\usepackage{amsmath}
\usepackage{tikz}
\usepackage{float}

\input amssym.def 
\input amssym.tex 

\newcommand{\howmany}{$28$ }

\begin{document}

\title{New Lower Bounds for 28 Classical Ramsey Numbers}

\author{
Geoffrey Exoo \\
\small Department of Mathematics and Computer Science \\[-0.6ex]
\small Indiana State University \\[-0.6ex]
\small Terre Haute, IN 47809 \\
\small\tt ge@cs.indstate.edu \\
\\
Milos Tatarevic \\
\small Alameda, CA 94501 \\
\small\tt milos.tatarevic@gmail.com \\
}

\date{}
\maketitle

\begin{center}
\end{center}

\begin{abstract}

We establish new lower bounds for \howmany classical two and three color
Ramsey numbers, and
describe the heuristic search procedures used.
Several of the new three color bounds are derived from the
two color constructions; specifically, we were able to use $(5,k)$-colorings
to obtain new $(3,3,k)$-colorings, and $(7,k)$-colorings to obtain new $(3,4,k)$-colorings.
Some of the other new constructions in the paper are derived from two
well known colorings: the Paley coloring of $K_{101}$
and the cubic coloring of $K_{127}$.

\end{abstract}

\section{Introduction}

The classical Ramsey number $R(k_{1},k_{2},\dots,k_{r})$ is the smallest integer $n$
such that in any $r$-coloring of the edges of the complete graph $K_n$ there
is a monochromatic copy of $K_{k_{c}}$ for some $1\leq c \leq r$.
In this note, we examine colorings for $r=2$ and $r=3$,
describe some new variations on computer construction methods
and are able
to improve the lower bounds for \howmany classical Ramsey numbers.
When discussing the procedures for improving lower
bounds, we use the term {\em bad subgraphs} to mean monochromatic complete graphs of
order $k_{c}$ in any of the colors $1\leq c \leq r$.
The reader is referred to Radziszowski's survey on
Small Ramsey Numbers \cite{survey} for basic terminology related
to the problem. The survey contains a comprehensive summary
of the current state of the art.

In this report the emphasis is on establishing new lower bounds for a range of Ramsey numbers
of the form $R(k_1,k_2)$.
Our primary focus is on cases where $4 \leq k_1 \leq 6$
and where current lower bounds are between $100$ and $200$,
though several bounds outside this range are also presented.
Some of these colorings are then used to obtain lower
bounds for three-color classical Ramsey numbers.
The new values for two- and three-color numbers are listed in Tables 1 and 2,
respectively.

The graphs were constructed by three variations on a basic method.
These three variations involve searching through colorings such that the
adjacency (coloring) matrix is partitioned into
(1) cyclic orbits, (2) square sub-blocks with cyclic orbits,
or (3) single edge orbits.
In the first case we get circle colorings, which have been used
many times in this context \cite{survey}.
In the second case we have colorings such that the adjacency
matrix can be partitioned into square circulant submatrices (of the same size);
and in the third case we have arbitrary coloring matrices.

Given a particular Ramsey number, $R(k_1,k_2)$ and a value of $n$, the goal is
to find a good $(k_1,k_2)$-coloring of $K_n$.
We begin by searching for circle colorings.  If that search succeeds,
then we might increase $n$, and repeat.  If it doesn't succeed and if $n$ is not
prime, then we expand our search space by considering all integer
factorizations $n = m d$, and search
for circulant block colorings where the adjacency matrix is partitioned into
an $m \times m$ array of $d \times d$ circulant blocks.

If neither of these searches succeeds, but if either of them produces colorings with a
relatively small number of bad subgraphs, then the coloring is used as a starting
configuration and a search which recolors single edges is applied.  The decision as to
whether a given number of bad subgraphs is small enough to merit further
investigation is made empirically.

In a further variation on this method,
a few of the new colorings presented here were obtained by using one of two well known
circulant colorings to create initial colorings, and then applying the single edge
recoloring procedure.
These two special colorings are the
Paley (or quadratic) coloring of $K_{101}$ and the cubic
coloring of $K_{127}$.
Recall that
the former coloring was used to establish the current lower bound for
$R(6,6)$ \cite{kalb}, while 
the latter was used to establish lower bounds for both
$R(4,12)$ \cite{sll} and $R(4,4,4)$ \cite{hillirving}.
The cubic residue graph of order $127$ may be of further interest, since
it was conjectured by the first author to be a Folkman graph
i.e., a $K_4$-free graph such that in any two-coloring of the edges, there
is a monochromatic $K_3$ (see, for example, Conjecture 4.4 in \cite{xurad}
which contains a fascinating discussion of this and related problems).

\section{New Lower Bounds}

As indicated above, we present two tables, one each for two-color and
three-color lower bounds.
A few of the cases listed were obtained using {\it ad hoc} methods,
which are outlined below.
For the other cases, we stayed fairly
close to the general approach described above.
In the column labeled {\it method} we indicate how we obtained the initial coloring
that gave, or led to, the indicated lower bound.
Since these were computer
searches, and the success of such searches depends to some extent on the
computer time available, we attempted to assign approximately the same
time to each problem.

\begin{table}[H]
\centering
\begin{tabular}{lccc}
Ramsey Number & Old Bound & New Bound & Method \\ \hline
R(4, 8) &  58 &  59 & circulant blocks \\
R(4,11) &  98 & 102 & circulant blocks \\
R(4,13) & 133 & 138 & cubic(127) \\ 
R(4,14) & 141 & 147 & cubic(127) \\ 
R(4,15) & 153 & 155 & cubic(127) \\
R(4,16) & 164 & 166 & cubic(127) \\ \hline
R(5,10) & 144 & 149 & circulant-minus  \\
R(5,11) & 171 & 174 & circulant  \\ 
R(5,12) & 191 & 194 & circulant  \\ 
R(5,13) & 213 & 218 & circulant  \\ 
R(5,14) & 239 & 242 & circulant  \\ 
R(5,15) & 265 & 269 & circulant  \\ 
R(5,16) & 290 & 293 & circulant  \\ \hline
R(6, 7) & 113 & 115 & Paley(101) \\
R(6, 8) & 132 & 134 & Paley(101) \\
R(6, 9) & 169 & 175 & circulant \\ 
R(6,10) & 179 & 185 & circulant \\ \hline
R(7, 9) & 241 & 242 & circulant  \\ \hline
\end{tabular}
\caption{New lower bounds for 2-color classical Ramsey numbers.}
\end{table}

In Table 2, new lower bounds for $3$-color Ramsey numbers are shown.
Some of these colorings were found by using one of our $(5,k)$-colorings,
and {\em splitting} the $K_5$-free color graph into two $K_3$-free color graphs,
therby obtaining a $(3,3,k)$-coloring.
In two other cases we used one of our
$(7,k)$-colorings to obtain a $(3,4,k)$-coloring.  Here we split
the $K_7$-free color graph into a $K_3$-free graph and a $K_4$-free graph.

\begin{table}[H]
\centering
\begin{tabular}{lccc} 
Ramsey Number & Old Bound & New Bound & Method \\ \hline
R(3,3,6) &   60 &   61 & circulant blocks \\
R(3,3,10) &  147 &  150 & circulant \\
R(3,3,11) &  162 &  174 & splitting \\
R(3,3,12) &  185 &  194 & splitting \\
R(3,3,13) &  212 &  217 & splitting \\
R(3,3,14) &  233 &  242 & splitting \\
R(3,3,15) &      &  269 & splitting \\
R(3,3,16) &      &  291 & splitting \\ \hline
R(3,4,7) & 145 & 152 & splitting \\
R(3,4,9) & 229 & 242 & splitting \\ \hline
\end{tabular}
\caption{New lower bounds for 3-color classical Ramsey numbers.}
\end{table}

All of the colorings described in this paper are available at
the Electronic Journal of Combinatorics and also at
the following location.

\begin{center}
\verb+http://cs.indstate.edu/ge/RAMSEY/ExTa+
\end{center}

\section{Search Methods}

Our basic method begins by considering colorings where
the adjacency matrix is either a circulant or can be partitioned into square
circulant submatrices of the same size.
For a graph of order $n$, let $m$ be a positive integer such that $m\mid n$.
We can partition the coloring matrix $A$ into $m \times m$ square blocks, as follows

\begin{equation}
\label{block}
A=\left[\begin{array}{cccc}
C_{1} & C_{2} & \cdots & C_{m}\\
C_{2}^{t} & C_{m+1} & \cdots & C_{2m-1}\\
\vdots & \vdots & \ddots & \vdots\\
C_{m}^{t} & C_{2m-1}^{t} & \cdots & C_{\binom{m+1}{2}}
\end{array}\right],
\end{equation}

\noindent and constrain each of the blocks $C_{i}$ to be a $d\times d$
circulant matrix, with $d=n/m$.
Note that the diagonal blocks must be symmetric circulant matrices.

For $m=1$, the entire coloring matrix of $K_{n}$ is a symmetric circulant. 
For $m>1$, the off-diagonal matrices can be either symmetric
or asymmetric circulant blocks.
The former case significantly reduces the size of the
search space, which makes things faster, but may eliminate some
good colorings.
When off-diagonal blocks are asymmetric, the number of required color
choices is

\begin{equation}
\label{block-sym}
m \lfloor d/2 \rfloor + \binom{m}{2} d.
\end{equation}

Whereas, if they are required to be
symmetric, the number of color choices is

\begin{equation}
\label{block-asym}
m \lfloor d/2 \rfloor + \binom{m}{2} (\lfloor d/2 \rfloor + 1).
\end{equation}

Let $V$ denote the vector of color assignments
that determines a coloring of $K_{n}$ where the
adjacency matrix is partitioned into circulant blocks, as in (\ref{block}).
When $m=1$ (and we have a circle coloring), the coloring can be 
completely specified by a vector $V=(a_{1},\dots,a_{\lfloor n/2 \rfloor})$.
For $m>1$,
the length of the coloring vector $V$ is either (\ref{block-sym}) or (\ref{block-asym}).
In the case where we consider recoloring single edges, the coloring vector
$V$ has length $\binom{n}{2}$.

In the algorithms described below, we use the notation $V_{i,c}$ to indicate
the vector obtained from
a coloring vector $V$ by assigning color $c$ to component $i$ of vector $V$.

For a given coloring vector $V$,
let $f_{c}(V)$ denote the number of monochromatic complete graphs of order $k_{c}$ in
color $c$.
During the search, the goal is to minimize $f_{c}$ for each $ 1 \leq c \leq r$;
but attempting to minimize a simple count of bad subgraphs does not appear to be
the best strategy.
So we adjust our objective function by defining
the $score$ of a coloring vector $V$ as a weighted sum

\[
score(V) = \sum_{c=1}^{r} w_{c} \cdot f_{c}(V)
\]

\noindent
where $w$ is a scalar array of {\em weights} adjusted
during the course of the algorithm.
The array $w$ is useful for two reasons.
The first reason applies to off-diagonal Ramsey numbers and deals with the
relative difficulty of eliminating badgraphs in the two colors.
Consider the case of an off diagonal Ramsey number, $R(s,t)$, where $s < t$.
One will find that it is
easier to eliminate the monochromatic $K_t$ subgraphs than to eliminate
the monochromatic $K_s$ subgraphs.
Local search procedures will tend to generate colorings with a large number of
bad $K_s$ subgraphs and few (or no) bad $K_t$ subgraphs.  This tendency
increases as $|s - t|$ increases.
To combat this problem, one can use weights such that $w_{s} > w_{t}$.
A good rule of thumb is to begin the search procedure with weights that satisfy
\[
 \frac{t}{s} \leq \frac{w_{s}}{w_{t}} \leq \Big(\frac{t}{s}\Big)^2.
\]

A second reason for using weights is that they help avoid getting
stuck in local minima.
As the search proceeds, we adjust the weights so they are
proportional to a moving average
of the bad subgraph counts.
For the problems in this paper, we updated
$w$ for each color $c$ by
\[
  w_c = \frac{K w_c + f_c}{(K+1)\sum_{i=1}^{r} f_i}
\]
\noindent
where $K$ is a positive constant.
The choice of $K$ seems important.
When $K$ is too small, the weights will oscillate wildly, and so will the
number of bad subgraphs.
On the other hand, if $K$ is too large, the weights
will change too slowly, and the algorithm is susceptible to getting trapped in local
minima.
Values of $K$ between $10$ and $100$ seem to be effective for the
problems considered in this paper.

In addition, the weights can be changed randomly.  The relative size of the
random contribution to the weights is analogous to the temperature variable
in simulated annealing \cite{metro}.

In the discussions below, we have occasion to mention the CPU
time required to obtain a certain bound.
All such times are approximate
times needed for a single CPU core,
running at $3.0$ gHz (an average).
Several bounds obtained here required months of CPU time.
In these cases, the work was performed
using $288$ CPU cores at the Department of Mathematics and Computer Science at the
Indiana State University.  The majority of these CPUs are third generation Intel {\em i5}'s.
With this setup, a month of CPU time can be covered in less than three hours.

\subsection{Method 1}

Our basic search method combines steepest descent with tabu search \cite{tabu}.
An outline of the search procedure is displayed in the Algorithm 1.  We
emphasize that the coloring vector $V$ can be any one of the three types
mentioned above: circulant, block-circulant, or individual edges.

The procedure begins with
a random coloring vector (which determines the adjacency coloring matrix)
and maintains a record of the last $L$ colorings. 
During each iteration, the procedure considers each of the possible recoloring vectors
$V_{i,c}$ not in the tabu list.
A recoloring vector that produces the minimum score is chosen and applied.

If $L$ is large relative to $n$,
it is possible for the algorithm to enter a state in which no recoloring
$V_{i,c}$ can be used
without revisiting a state on the tabu list.
If this (rare) event occurs, we restart with a new random
coloring.  When coloring graphs of order $n$, where $150 < n < 200$
(which accounts for several of the cases considered here), setting $L = 1000$
seemed to produce good results.

\vspace{3mm}
\begin{algorithm}[H]
\caption{Steepest descent with tabu search}
\vspace{1mm}
\KwIn{Coloring vector $V$,
number of colors $r$,
size of a tabu list $L$}
\Begin{
  Set elements of $V$ to random colors\;
  $H \gets \{\}$\;
  \Repeat{$ score(V)=0 $  {\bf or} $ M=\emptyset $} {
    $W \gets \{\}$\;
    \For{$1 \leq i \leq \left\vert V \right\vert $}{
      \For{$1 \leq c \leq r$}{
        \If{$V_{i, c} \not\in H$} {
          $W = W \cup V_{i, c}$\;
          Compute $ score(V_{i, c}) $\;
        }
      }
    }
    $M \gets \{V \in W \,|\, score(V) \mbox{ is a minimum}\}$\;
    \If{ $ M \ne \emptyset $ } {
      Randomly choose new $V \in M$\;    
      $H \gets H \cup V$\;
      \If{$\left\vert H \right\vert > L$} {
        Remove an element from $H$ in FILO manner\;
      }
      Adjust weights\;
    }
  }
  \Return{$V$}\;
}
\end{algorithm}

\subsection{Method 2}

In addition to the weight based
method given in Algorithm 1,
we also used a procedure based on simulated annealing \cite{metro}.
This method is not suitable for circulants, but performed well
for several block-circulant colorings, and is a good alternative
when the recoloring of the individual edges is considered.
This local improvement search procedure
is outlined in Algorithm 2. At the beginning of the procedure,
we find the set of all edges contained in bad subgraphs
of order $k_{c}$ in any of the colors $1\leq c \leq r$.
This edge set is denoted by $I$ and can be created
during the bad subgraph counting procedure, using
very little additional CPU time. 
In the inner loop of Algorithm 2,
we apply simulated annealing on edges from the set $I$.
The initial temperature $T_{0}$ and the cooling rate depends
on $n$ (the order of the graph) and the initial $score$.
The inner loop is limited by $j_{max}$ iterations, after which we update the set $I$.
Typically we would set $j_{max} = \left\vert V \right\vert / \, 4$.
In most cases when $T$ reached $0$ the coloring can be improved further, so
we restarted the procedure using the newly obtained coloring vector $V$. The process is
repeated until the further improvements can not be made.

To perform a circulant block search, several modification are required.
Initially $V$ is set to a random coloring vector.
Next, we set $I = \{1,\dots,\left\vert V \right\vert\}$.
The optimal choice of the initial temperature $T_{0}$
and cooling rate varies greatly for different $n$, $m$ and $r$,
and we attempt to determine them empirically.
In practice, if $score(V)$ is decreasing more slowly than expected,
we restart the entire procedure using the same initial parameters.

\vspace{3mm}
\begin{algorithm}[H]
\caption{Search procedure based on simulated annealing}
\vspace{1mm}
\KwIn{Coloring vector $V$,
number of colors $r$,
initial temperature $T_{0}$,
number of iterations of the inner loop $j_{max}$}
\Begin{  
  $T \gets T_{0}$\;
  $s\gets \infty$\;
  \Repeat{$ s=0 $ {\bf or} $T = 0$} {
    $I \gets \{\}$\;
    \For{$1 \leq i \leq \left\vert V \right\vert $}{
      $c \gets V[i]$\;
      \If{$f_c(V_{i,c+1}) < f_c(V)$}{
        $I \gets I \cup \{i\}$\;         
      }
    }
    \For{$1 \leq j \leq j_{max} $}{
      Randomly choose $i \in I$\;	
      Randomly choose $c \in \{1,\dots,r\}$, $c \ne V[i]$\;
      $s^{*} \gets score(V_{i,c})$\;
      \If{$s^{*} - s \leq 0$ {\upshape{\textbf{or}}} $e^{(s - s^{*})/T} > random(0,1)$}{
        $V \gets V_{i,c}$\;
        $s \gets s^{*}$\;
      }
    }
    Reduce $T$\;
    Adjust weights\;
  }
  \Return{$V$}\;
}
\end{algorithm}

\subsection{Special Constructions}

Several lower bounds presented in this paper required a combination of our basic search methods.
A brief summary of how these colorings were obtained is given below.

\subsubsection{$R(5,10)$}

The lower bound for $R(5,10)$ (labeled {\it circulant-minus} in Table 1) required two steps.
First, we began with a circulant coloring of $K_{149}$ that was nearly a good coloring.
This particular coloring had $149$
monochromatic $K_5$'s and no monochromatic $K_{10}$'s.
After we had applied a local search procedure, one problem vertex was deleted,
and a local search for a good coloring succeeded.

\subsubsection{$R(4,8)$}

To improve the lower bound on $R(4,8)$, we examined colorings of $K_n$ for $60 \leq n \leq 64$,
and for all $m$ such that $m\mid n$ and $2 \leq m \leq 10$.
The circulant block search was performed using the simulated annealing procedure.
We then performed local searchs on those graphs we found with fewer than $250$ bad subgraphs.
For $m = 4$ and $m = 5$, graphs on $60$ vertices were found and used
to extract a $(4,8)$-coloring on $58$ vertices.
In the first case, the block construction contained $60$ monochromatic $K_4$'s
and $45$ monochromatic $K_{8}$'s.
In the second, it contained $63$ monochromatic $K_4$'s and no monochromatic
$K_{8}$'s.
In both cases, the best coloring was reached after removing two vertices.
We note that these graphs on $60$ vertices were not those with the minimum
possible number of bad subgraphs (among block circulant colorings).
Many colorings were found with fewer than $100$ bad subgraphs,
and some with fewer than $50$, but we were not be able to modify these
to obtain good colorings.

To obtain a new $(4,8)$-coloring on $K_{58}$ vertices, several hours of CPU
time were required. We spent over a month of CPU time trying to improve this solution.
We found several close graphs on $K_{59}$ and $K_{60}$. On $K_{59}$, we
were able to find a coloring with only $3$ monochromatic $K_4$'s and no monochromatic
$K_{8}$'s.

\subsubsection{$R(4,11)$}

Perhaps our most complicated search was for the case of
$R(4,11)$.
This search required four steps, beginning
with a $(4,8)$-coloring, extending it to a $(4,9)$-coloring, then to a
$(4,10)$-coloring, and finally to a $(4,11)$-coloring.
We started from a circulant block construction for $n=60$ and $m=4$,
the same one we had used for $R(4,8)$.
We extended this graph by appending circulant $15\!\times\!15$-blocks to
the $60\!\times\!60$ coloring, and without changing
the coloring on the first $60$ vertices.
A promising $(4,9)$-coloring on $75$ vertices was found.
Next we added another layer of $15\!\times\!15$-blocks and obtained an even more promising
$(4,10)$-coloring with $n=90$ and $m=6$.
Then we extended this coloring to a (nearly good) $(4,11)$-coloring on
$105$ vertices.  After repeatedly reaching local minima, we became convinced
that we were not going to find a good coloring 
while maintaining the circulant block structure.
At that point we removed four vertices and applied a local search procedure
recoloring individual edges, and obtained a good coloring on $101$ vertices,
thus establishing the new lower bound for $R(4,11)$.
The promising coloring on $105$ vertices required
several hours of CPU time to locate.
The final local search required
approximately a week of CPU time.

\subsubsection{$R(3,3,6)$}

A graph that we find interesting was uncovered when applying
a circulant block search to the three color Ramsey number $R(3,3,6)$.
The best result was obtained for $n=60$ and $m=6$, and no
additional corrections were required.
The coloring has significant symmetry.
The graph in the first color has $2^{23} \cdot 5 = 41943040$ automorphisms;
the color two graph $80$ automorphisms;
and the color three graph has $40$ automorphisms.

This coloring was discovered using the procedure similar to the one described in Algorithm 1,
and the search took only a few minutes of CPU time.
It is interesting that several other search procedures we
tried were not able to find a good coloring on $K_{60}$.

\subsection{Using Colorings of Lower Order}

In certain cases, when two Ramsey numbers $R(s,t)$ and $R(s,t+k)$ are close,
we tried to obtain a new lower bound for $R(s,t + k)$ by using the coloring that established
the lower bound for $R(s,t)$ as a starting point for the $(s,t+k)$ search.
Copies of $(s,t)$-coloring were used as starting blocks in
a block circulant coloring that we hoped to manipulate into
an $(s,t+k)$-coloring.
If the order of the desired coloring was not an integer
multiple of the order of the good $(s,t)$ coloring,
we deleted rows and columns until we had a matrix of the desired size.
If the resulting graph had only a few bad subgraphs,
we tried to improve it either by using a
block circulant search or a single edge recoloring procedure.

\subsubsection{Colorings derived from Paley coloring of $K_{101}$}

To obtain the new bounds for $R(6,7)$ and $R(6,8)$ we used the Paley coloring of $K_{101}$,
which was the graph used to establish the longstanding lower bound for $R(6,6)$ \cite{kalb}.
In both cases, the initial block coloring created from $101 \times 101$ circulant blocks was
reduced to a matrix of the desired size, and the resulting coloring contained a relatively large
high number of bad subgraphs.
By applying the block circulant search procedure, we reduced the
number of bad subgraphs significantly.
To obtain good colorings, several days of CPU time
were required. Additionally, several months of CPU time were spent trying to
further improve the bounds.
During this final stage, both bounds were improved by one.

Let us note that a similar result can be obtained for $R(6,7)$ by using only
the single edge recoloring procedure,
but the search required more CPU time.
For $R(6,8)$, it was necessary to make improvements in the block circulant structure first,
and then apply single edge recoloring.

\subsubsection{Colorings derived from a cubic coloring of $K_{127}$}

The new lower bounds for $R(4,t)$ for $13\leq t \leq 16$ were obtained by using the
cubic coloring of $K_{127}$. The cubic coloring was used to establish the lower bound
for R(4,12) \cite{sll}.
The previous bound for $R(4,13)$ was only $5$ greater than that for $R(4,12)$,
and in the case of $R(4,14)$ the difference
was only $13$.
It was expected that the search would succeed at least for $t=13$, but 
surprisingly, we quickly obtained the other good colorings.

The initial block coloring created from copies of the
cubic residue graph on $K_{127}$ contained only few bad subgraphs. 
However, performing a local search on those colorings was particularly slow.
In all four cases, we feel that further improvement is possible.
But better bounds would likely require years of CPU time using this method.

\subsection{Splitting}

This section describes a simple method that allowed us to reach most of the lower
bounds given in Table 2.
The idea is to use known $(R(s,t)-u,k)$-colorings to create $(s,t,k)$-colorings, for some
$u\ge 1$.
We noticed that some of the bounds for Ramsey numbers of the form $R(3,3,k)$
were close to the corresponding bounds for numbers of the form $R(5,k)$.
A similar observation was made with regard to $R(3,4,k)$ and $R(7,k)$.  
In the first case, $5$ is just one less than the Ramsey number $R(3,3)$, and in the
second case $7$ is two less than $R(3,4)$.
So it appeared conceivable that we could use
known $(5,k)$-colorings to obtain new $(3,3,k)$-colorings, and similarly
use $(7,k)$-colorings to obtain new $(3,4,k)$-colorings.
In the first case, we have to split the first color into two $K_3$-free graphs.
In the second case, we have to split the first color into a $K_3$-free
graph and a $K_4$-free graph.

As a starting point we used $2$-color circulant colorings
obtained by the method described in Section 3.1.
First, we attempt to split the color one graph into two circulants.
If good coloring can not be achieved,
we continue the process by recoloring individual edges.
At this point, we still do not modify the third color.
In most of the cases, we found good colorings almost instantly.

Of course there is no guarantee that this procedure will work in general.
In some cases, we were not able to find a good coloring using a given circulant, so
we generated other non-isomorphic circulant colorings and tried each of them.
Usually after few attempts the procedure succeeded.
In several cases, we were unable to find good circulant colorings.
But if the number of bad subgraphs was small enough,
we ran the single
edge recoloring procedure again, this time allowing changes to the
third color.
We were thereby able to find,
for example, a good $(3,3,12)$-coloring on $193$ vertices.
But the process took a few hours of CPU time.
Note that the lower bounds we found for $R(3,3,13)$ and $R(3,3,16)$
are slightly less than those found for $R(5,13)$ and $R(5,16)$, reflective of the
fact that we could not complete the splitting procedure for our best $(5,13)$ and
$(5,16)$ colorings, but did succeed using smaller colorings.

\section{Conclusions}

One of the motivations of our work was to apply a fixed amount resources,
both computer and human, to a range of Ramsey problems, and learn to what
extent the lower bounds could be improved.  We focused on classical two-color
Ramsey numbers where the best known lower bound was between $100$ and $200$.
We did this for two reasons.  First, neither of us had ever applied serious
effort to problems in this range, and second, we felt these would be problems
where progress could be made.  For problems where the current lower bound is
less than $100$, we felt that it was unlikely that progress would be made
by simple methods (the case of $R(4,8)$ turned out to be an exception).
For problems where the current lower
bound is significantly greater than $200$, we felt that the required computer
time would usually be too great.

In each of the cases considered, we applied all of the methods described here,
and few others.  Other methods tried, but not mentioned here include the use
of general (i.e., non-cyclic) Cayley colorings.   It is
remarkable that for every problem considered here, circulant colorings
do better (produce good colorings for a larger value of $n$) than Cayley
colorings from other groups.
One obvious explanation for this is that cyclic groups have elements of
larger average order than other groups, but it is still somewhat surprising.

Finally we considered the case of $R(8,8)$, where the current lower bound
seems much too small, as one may observe from the Figure (the data is taken from
\cite{survey}).  Despite signifcant computer time spent searching for
circulant colorings in the $282$ to $286$ range, we did not find anything
new.  It is tempting to conclude that if there is a circulant coloring that
improves the lower bound for $R(8,8)$, it is quite a bit larger than the
current best coloring.
In fact, even the related case of $R(7,8)$ is surprisingly difficult,
especially considering that the current bound of $217$ is only $12$ more
than the $R(7,7)$ bound.  We applied a variety of methods to this problem
with no success.

\vspace{5mm}

\begin{figure}[H]
\centering
\begin{tikzpicture}[scale=0.7]
\draw (2,1) -- coordinate (x axis mid) (11,1);
\draw (2,1) -- coordinate (y axis mid) (2,10.2);
\foreach \x in {3,...,10}
  \draw (\x,1.05) -- (\x,0.9)
    node[anchor=north] {\x};
\foreach \y in {2,...,10}
  \draw (2.05,\y) -- (1.9,\y)
    node[anchor=east] {\y};
\node[below=0.75cm] at (x axis mid) {n};
\node[rotate=90, above=0.75cm] at (y axis mid) {$log_2 R(n,n)$ \, (lower bound)};

\draw plot[mark=*,mark options={fill=blue}]
  file{log.data};

\draw[->,thick,red,>=stealth]
  (9,6) node[sloped,above,right] {R(8,8)} -- (8.06,8);

\end{tikzpicture}
\caption{Log plot of current lower bounds for $R(n,n)$}
\end{figure}
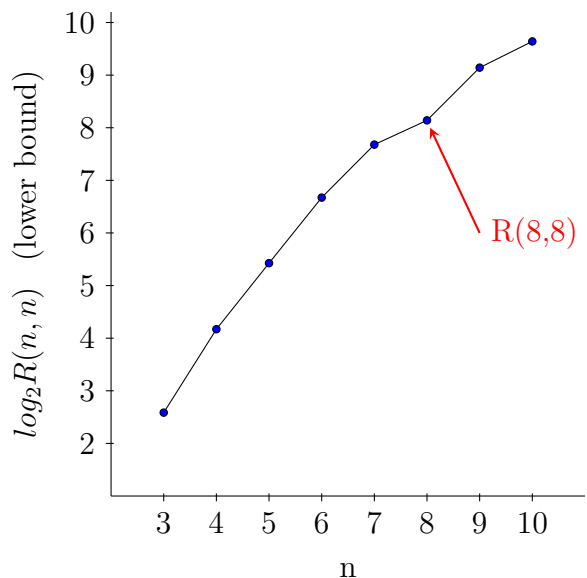

\bibliographystyle{plain}

\begin{thebibliography}{99}

\bibitem{tabu}
F. Glover,
\newblock Tabu Search -- Part 1,
\newblock {\em ORSA Journal on Computing}, 1: 190--206, 1990.

\bibitem{HaKr}
H. Harborth and S. Krause,
\newblock Ramsey Numbers for Circulant Colorings,
\newblock {\em Congressus Numerantium}, 161: 139--150, 2003.

\bibitem{hillirving}
R. Hill and R.W. Irving,
\newblock On Group Partitions Associated with Lower Bounds for Symmetric Ramsey Numbers,
\newblock {\em European Journal of Combinatorics}, 3: 35--50, 1982.

\bibitem{kalb}
J.G. Kalbfleisch,
\newblock Construction of Special Edge-Chromatic Graphs,
\newblock {\em Canadian Mathematical Bulletin} 8: 575--584, 1965.

\bibitem{metro}
S. Kirkpatrick, C. D. Gelatt and M. P. Vecchi,
\newblock Optimization by Simulated Annealing,
\newblock {\em Science}, 220: 671--680, 1983.

\bibitem{survey}
S. P. Radziszowski,
\newblock Small Ramsey Numbers,
\newblock {\em The Electronic Journal of Combinatorics}, DS1, 2014

\bibitem{xurad}
X. Xu and S. P. Radziszowski,
\newblock On Some Open Questions for Ramsey and Folkman Numbers,
\newblock {\em The Royal Swedish Academy of Sciences}, Report 19,
2013/2014.

\bibitem{sll}
S. Wenlong, L. Haipeng and L. Qiao,
\newblock New Lower Bounds of Classical Ramsey Numbers R(4,12), R(5,11) and R(5,12),
\newblock {\em Chinese Science Bulletin}, 43(6): 528, 1998.

\end{thebibliography}

\end{document}